\documentclass[12pt]{article}
\usepackage{MnSymbol}
\usepackage{amsmath,amsfonts,
verbatim}

\newcommand{\Q}{{\mathbb Q}}

\newcommand{\Z}{{\mathbb Z}}

\newcommand{\LL}{\mathcal{L}}

\newcommand{\F}{\mathbb{F}}

\newcommand{\X}{{\mathbb X}}

\newcommand{\pr}{{\rm pr}}

\newcommand{\M}{\mathbf{M}}   

\newcommand{\NN}{\mathbf{N}}

\newtheorem{pkt}{}[section]  
\newcommand{\bpk}{\begin{pkt}\rm }  
\newcommand{\epk}{\end{pkt}} 
\newcommand{\inv}{^{-1}}   

\newcommand{\be}{\begin{equation}}  
\newcommand{\ee}{\end{equation}}

\newcommand{\dcl}{\mathrm{dcl}}
\newcommand{\acl}{{\rm acl}}

\newcommand{\kk}{\mathrm{k}}

\newcommand{\Aut}{\mathrm{Aut}}
\newcommand{\tp}{\mathrm{tp}}
\usepackage{graphicx}

\usepackage{tikz}

\title{A model theory section conjecture}
\author{R.Abdolahzadi and B.Zilber\footnote{Supported by the EPSRC program grant ``Symmetries and Correspondences''}}
\date{ November 9, 2020}

\begin{document}
\maketitle

We  introduce the category of structures and interpretations which allows us to discuss some issues of Grothendieck's anabelian geometry in model-theory terms.

Our main result is a formulation in terms of pure stability theory of a problem closely related to Grothendieck's section conjecture.\footnote{We would like to thank J.Derakshan for his interest and many useful comments on the paper. }

\section{The category of strutures and interpretations}\label{section 1}
\bpk Most of the material below is known. See \cite{CDM} and \cite{MTB} for a model-theoretic approach which we further pursue here. The community of anabelian geometers prefers to speak in terms of Galois categories, see e.g. \cite{C}. One of the aims of the current project is to demonstrate advantages of the model-theoretic point of view. 

Unlike the above publications we do not apriori  restrict the power of the language to first order. The default assumptions is that the expressive power of the language is such that the following holds:

\be \label{def=inv} \begin{array}{ll}\mbox{\em A subset of a structure $\M$ is definable iff} \\
\mbox{\em it is invariant under  automorphisms of $\M$}
\end{array}
\ee

For finite structures this property holds for first order languages. For countable structures the language $L_{\omega_1,\omega}$ serves the purpose (by Scott's Theorem, \cite{Scott}). More generally we require a {\em powerful enough language} the choice of which depends on $\M.$ 

 The main interest to us are {\em finitary} structures defined below. For this class of structures first-order languages are essentially sufficient.

In general, we consider multi-sorted $\LL_\M$-structures $\M$. A {\bf definable set} in 
$\M$ is a definable (without parameters) subset $D$ of $\prod_{i\in I} M_i,$ some $I,$ a cartesian product  of sorts. When $\M$ stands for a structure, $M$ is the universe of the structure, the union of its sorts.
 
A definable {\bf sort} in an $\LL_\M$-structure $\M$ is a set of the form $D/E$ where $D=D(\M)$ is a definable set in $\M$ and $E$ a definable equivalence relation on $D.$ An $n$-ary relation  on  $D/E$ is definable if its pull-back under the canonical map
$D\to D/E$ is  definable. 

{\bf An interpretation} of an  $\LL_{\NN}$-structure $\NN$ in an $\LL_\M$-structure $\M$ is a bijection
$g: N\to D/E,$  a sort in $\M$ such that 
for each basic relation $R$ the image $g(R)$ is a definable relation on the sort  $D/E.$

Given a structure $\M$ we also consider the structure $\M^{Eq}$ 
(non-first-order version of $\M^{eq}$)
 interpretable in $\M$ and which has every sort of $\M$ as a definable substructure.

Note that any union of sorts and  a direct product of any number of sorts is a sort in $\M^{Eq}.$

We reserve the notation $\M^{eq}$ for  the extension of $\M$ by first-order imaginary sorts, see \cite{TZ}.
 
\epk
\bpk\label{facts} Standard facts about first-order imaginaries (see e.g. \cite{TZ}) easily generalise to $\M^{Eq}$ with the help of (\ref{def=inv}). 

To every relation $R$ definable in $\M$ using parameters is associated a {\bf canonical parameter} $c\in \M^{Eq}$ which is fixed by the same automorphisms  as fix $R.$ More generally, if $\NN$ is interpretable in $\M$ using parameters there is a canonical parameter for $\NN$ in $\M^{Eq}$
 which is fixed by exactly the automorphism of $\M^{Eq}$ which act on $\NN$ as automorphisms of $\NN.$ Canonical parameters are defined uniquely up to interdefinability, that is the condition  $\dcl(c)=\dcl(c').$ 
 
   We use 
 $$^\lceil{\NN}^\rceil=\{ c\in \M^{Eq}: \forall \sigma\in \Aut\, \M^{Eq}\, \sigma(c)=c \leftrightarrow \sigma_{|\NN}\in \Aut\, \NN\}$$
 to define the set of all interdefinable canonical  parameters.

 Note that in a powerful enough language (e.g. a second order language) the sorts
$$\mathcal{R}_n:=\{ R: R\subseteq M^n\}$$
are interpretable in $\M.$

In particular, an arbitrary subset or  relation $R$ on $M$ is interpretable in $\M$ using a parameter in a $\mathcal{R}_n.$ The same is true for structures.

\epk
\bpk {\bf Category $\mathfrak{M} .$} Its objects are (multisorted) $\LL_\M$-structures $\M$.
 
 The {\bf pre-morphisms} $g:\mathbf{N}\to \M$
 are interpretations (without parameters).  More precisely,
 $$g: \mathbf{N}\to \M^{Eq}$$ is an injective map such that $g N$  is a universe of a sort in $\M^{Eq},$  and for  any basic relation or operation $R$ on  $\mathbf{N}$ the image $g R$ is definable  in the sort. 
 
 We denote $g\mathbf{N}$ the  $g N$ together with all the relations and operations $g R$ for $R$ on $\NN.$

\medskip
 
Two pre-morphisms $g_1: \NN\to \M$ and $g_2\to \NN\to \M$ are equivalent if  there is a bijection $h:g_1\NN\to g_2\NN$   which is definable in $\M^{Eq}.$

The equivalence class of a pre-morphism $g:\NN\to \M$ is a {\bf morphism} $\mathbf{g}:\NN\to \M.$  

\medskip

The following definitions will be used  for pre-morphisms $g$  as well as for morphisms
$\mathbf{g}.$
 
 We say $g$ is an {\bf embedding},  $g:\mathbf{N}\hookrightarrow \M$ if $g\mathbf{N}$ has no proper expansion definable in $\M^{Eq}.$

We say $g$ is a {\bf surjection}, $g:\NN\twoheadrightarrow \M$ if $M\subseteq \dcl(g N)$
where $\dcl$ is in the sense of $\M^{Eq}.$

We say that $g: \NN\to \M$ is an {\bf isomorphism} and write  $g:\mathbf{N}\cong \M,$ if
$g$ is an embedding and a surjection.
 
\medskip

In what follows we sometimes write $\mathbf{N}\cong_{\mathfrak{M}} \M$ to emphasise that the isomorphism (or morphism) is in the sense of the category $\mathfrak{M}$ to distinguish from ones in the usual algebraic sense, and call it an $\mathfrak{M}$-isomorphism.

\epk
\bpk
{\bf Lemma.} {\em Let $g:\mathbf{N}\to \M$ be an $\mathfrak{M}$-isomorphism and let $\M'=g\NN.$ Then the inverse map $g^{-1}:\M'\to \mathbf{N}$ induces 
 a $\mathfrak{M}$-isomorphism 
$h: \M\to \NN.$}

\medskip

{\bf Proof.} 
By assumptions
we have  $M\subseteq \dcl (M')$ in $\M^{Eq}.$ This implies that there are in $\M:$
a family $\{ S_i: i\in I\}$   of definable
subsets $S_i\subset  {M'}^{n_i}$ and 
 a family of definable functions $h_i: S_i\to M$ such that $$\bigcup_{i\in I} h_i(S_i)=M\mbox{ and }h_i(S_i)\cap h_j(S_j)=\emptyset \mbox{ if } i\neq j .$$

{\bf Claim.} We may assume that the family $\{ S_i: i\in I\}$ of domains of $h_i$ is disjoint, that is $S_i\cap S_j=\emptyset$ if $h_i\neq h_j.$

Proof. Note that by definition $\dcl(M')=\dcl(M'\cup \dcl(\emptyset)),$ where $\dcl$ is understood in the sense of $\M^{Eq}.$ The latter has, for each $i\in I,$ the sort $^\lceil{h_i}^\rceil$ which is defined as the $\mathrm{graph} h_i/E_i$ where $E_i$ is the trivial equivalence relation with one equivalence class. Clearly, $^\lceil{h_i}^\rceil \in \dcl(\emptyset).$ 
 Now replace 
$S_i$ by $S_i\times ^\lceil{h_i}^\rceil$ and we have the required.
   
   \medskip

Set $D(\M'):= \bigcup_{i\in I} S_i$ and 
$h: D(\M')\to \M$ to be  $\bigcup_{i\in I}h_i.$
 $h$ is a map definable in $\M$ and is an interpretation, a pre-morphism $\M\to \M'.$ On the other hand, any relation  on $\M'$ is a relation on a sort in $\M$ since $\M'$ is a sort in $\M^{Eq},$   hence there are no new relations on $h(M)$ induced from $\M',$ that is the interpretation $h$ is an embedding.  Recalling that $\M'= g \NN$ completes the proof. 
 $\Box$ 

\medskip

We identify morphism $h$ as in the Lemma with $g\inv.$
 \epk
\bpk For a subset $A\subseteq M,$ denote $\M/A$  the expansion of $\M$ by  names of elements of $A.$

Clearly, the identity map defines a (canonical) morphism $\M\to \M/A.$ This morphism is an embedding (and so isomorphism) if and only if $A\subseteq \dcl(\emptyset).$
\epk
\bpk \label{dcl} Given $A\subseteq \dcl(\emptyset)$ we may treat $A$ as a structure in which any element is named (e.g. by a formula defining the element in $\M$) and so any relation is definable. Clearly then
$$\mathrm{Aut}(A)=1 \mbox{ and } A\hookrightarrow \M.$$
 
\epk
\bpk {\bf A remark on notation.} The category $\mathfrak{M}$ treats
$\M$ and $\M^{Eq}$ as isomorphic objects, so we often do not distinguish between the two in our notation. In this context the notation $\M/A$ makes sense even when $A\subset \M^{Eq}.$

\epk 
 
 \bpk\label{finry} The category  $\mathfrak{M}_\mathrm{fin}$ is a subcategory of  $\mathfrak{M} $ whose objects are {\bf finitary} structures $\M,$ that is structures which can be represented in the form
 $$\M=\bigcup_{\alpha<\kappa } \M_\alpha$$
 where the $\M_\alpha$ are finite first-order $0$-definable substructures of $\M.$
 \medskip

 Note that an equivalent definition would be 
 $$\M=\acl(\emptyset)$$
 where $\acl$ is in the sense of first-order logic.
\medskip

 {\bf Example.} Let $\kk$ be a field and $\F=\tilde{\kk},$ its algebraic closure. We consider $\F=\F/\kk$ as a structure
 in the language of rings with names for elements of $\kk.$ Then each $a\in \F$ is contained in a $0$-definable set $M_a$ equal to its Galois orbit $M_a:=G_k\cdot a,$ $G_k=\mathrm{Gal}(\F:\kk).$  So  $\F/\kk\in \mathfrak{M}_\mathrm{fin}.$
 \epk
 
 \bpk \label{EqCat}{\bf Theorem.} {\em The map $\M\mapsto \mathrm{Aut}(\M)$ induces a contravariant  functor from
 $\mathfrak{M}$ into the category  $\mathfrak{G}$ of  groups.  
 
(i) To every $g:\mathbf{N}\to \M$ corresponds the restriction homomorphism 
 
\noindent $\hat{g}:\mathrm{Aut}(\M^{Eq})\to \mathrm{Aut}(\mathbf{N}).$
 
 (ii) An embedding $g:\mathbf{N}\hookrightarrow \M$ corresponds to the surjection $\hat{g}:\mathrm{Aut}(\M)\twoheadrightarrow \mathrm{Aut}(\mathbf{N}).$ 
 
 (iii) The expansion by naming all points in $A\subseteq \M^{Eq},$
 $g:\M\to \M/A$ corresponds to an embedding $\hat{g}:\mathrm{Aut}(\M/A)\hookrightarrow \mathrm{Aut}(\M).$

(iv) The restriction of  $\Aut$   to the finitary subcategory $\mathfrak{M}_\mathrm{fin}$ sends into the category of profinite topological groups and continuous homomorphisms $\mathfrak{G}_{\mathrm{pro}}.$ Moreover, 
$$\mathrm{Aut}_\mathrm{fin}:\ \mathfrak{M}_\mathrm{fin}\to \mathfrak{G}_{\mathrm{pro}},$$ 
 is an equivalence of categories. 
 

   
}
\medskip

{\bf Proof.} (i)
is immediate by definition.

(ii) Since $g$ is an embedding, the relations definable on $g\mathrm{N}$ are the same in $\M$ and $\mathbf{N}.$ Hence a $g\mathbf{N}$- automorphism $\alpha$ is a
monomorphism (in the sense of a powerful enough language) $g\mathrm{N}\to g\mathrm{N}$ in $\M.$ Now use the transfinite back-and-forth induction  to extend $\alpha^*$ to a monomorphism $gN\cup M\to gN\cup M,$ equivalently, an automorphism of $\M.$ Clearly, $\hat{g}\mapsto \alpha.$
 

(iii) Immediate.

(iv) 
First we prove the statement for $\Aut:\mathfrak{M}_\mathrm{finite}\to \mathfrak{G}_{\mathrm{finite}},$ the functor between finite structures and finite groups, 
subcategories of  $\mathfrak{M}_\mathrm{fin}$ and  $\mathfrak{G}_{\mathrm{pro}},$ respectively.

Given a finite group $\mathbf{G}$ one constructs a finite $\M$ such that $\mathbf{G}\cong \mathrm{Aut}(\M)$ by setting $M=G$ and introducing all relations $R$ on $M$ which are invariant under the action of $\mathbf{G}$ on $G$ by multiplication. This gives us  $\M=(M; R)$ 

Claim.
$$\mathbf{G}=\mathrm{Aut}(\M)$$
Proof. $\mathbf{G}$ acts on $\M$ by automorphisms by definition. We need to prove the inverse, i.e. that there are no other automorphisms. Consider the tuple $\bar{g}$ of all the elements of $G$ (of length $n=|G|$) and let $S_g$ be the conjunction of all the relation in $R$ that hold on $\bar{g}$ (that is $\mathrm{tp}(\bar{g})$). We can also consider $S_g^0:=\mathbf{G}\cdot \bar{g},$ the orbit of $\bar{g}$ under the action of $\mathbf{G}.$ Clearly,
$S_g^0\subseteq S_g$ and by minimality they are equal. 

Now take an automorphism $\sigma$ and consider $\sigma \bar{g}.$ This is in $S_g$ and thus, for some $h\in G,$    $\sigma \bar{g}=h\bar{g},$ that is $\sigma g_i=hg_i$ for each $g_i\in G.$
Claim proved.

It remains to see that if  $\mathbf{G}\cong\mathrm{Aut}(\mathbf{N}),$ then   $\mathbf{N}$ is definable in $\M$ and vice versa. In order to do this we may assume $\mathbf{G}=\mathrm{Aut}(\mathbf{N}).$  

Consider $N,$ the universe of the structure, and let $\mathbf{n}$ be the $N$ presented as an ordered tuple. 
Let $M':= \mathbf{G}\cdot \mathbf{n},$ the orbit of the tuple under the action of the automorphism group. Clearly, $M'$ consists of $|\mathbf{G}|$ distinct elements, since automorphisms differ if and only if they act differently on the domain $N.$ Also $M'$ is definable in $\mathbf{N}$ since the tuples $\mathbf{n}'$ making up $M'$ are characterised by the condition that $\mathrm{tp}(\mathbf{n}')=\mathrm{tp}(\mathbf{n}).$ The relations $R$ induced on $M'$ from $\mathbf{N}$ are invariant under $\mathrm{Aut}(\mathbf{N}),$ and because a finite structure is homogeneous, the converse holds. In other words an obvious bijection $M\to M'$ is a bi-interpretation, so $\M\cong \M'$ in the sense of $\mathfrak{M}.$  At last notice that we can interpret $\mathbf{N}$ in $\M'$ since the relation  ``$\mathbf{n}'$ and $\mathbf{n}''$ have the same first coordinate is invariant under $\mathbf{G}$'' is definable. This gives us $N$ as a definable sort. It follows that any relation on $N$ definable in $\mathbf{N}$ is definable in $\M'.$ So $\mathbf{N}\cong \M'\cong \M$  in the sense of $\mathfrak{M}.$  Finite case of $\Aut$ proven.

Now we extend $\Aut$ to the category of  finitary structures   $\M\in \mathfrak{M}_{\mathrm{fin}}$ by continuity
$$\M=\lim_\to \M_\alpha \to \mathbf{G}=\lim_{\leftarrow}  \mathbf{G}_\alpha, \mbox{ where } \mathbf{G}_\alpha=\Aut  \M_\alpha$$

Since the functor is invertible and preserves morphisms on finite objects of the categories, it is an equivalence also on the limits.

$\Box$
 \epk
 \bpk {\bf Example.} Let $K, L\subset \bar{\Q}$ be two number fields and
   $\bar{\Q}_K$ and $\bar{\Q}_L$ be two structures, the fields of algebraic nubers with respective subfields of constants (named points). Clearly 
these belong to $\mathfrak{M}_\mathrm{fin}.$
A celebrated theorem by  Neukirch/Uchida \cite{Uchida} states that
$$\mathrm{Gal}(\bar{\Q}: K)\cong \mathrm{Gal}(\bar{\Q}: K)\Leftrightarrow K\cong L $$
and hense by \ref{EqCat}(iv)
$$\bar{\Q}_K\cong_\mathfrak{M}\bar{\Q}_L\Leftrightarrow K\cong L.$$

\epk

 \bpk\label{M/A} {\bf Lemma.} {\em Suppose $\NN,\M\in \mathfrak{M}$ and
  $$\hat{g}:\Aut\, \NN\hookrightarrow \Aut\, \M.$$
  Then there is $A\subset \M^{Eq}$ such that
  $$\NN\cong_\mathfrak{M} \M/A.$$}
  
  {\bf Proof.}
 By \ref{EqCat}(iv) we have $$g: \M\twoheadrightarrow \NN,$$ that is
 $g\M$ is a substructure of $\NN^{Eq}$ such that $\dcl(gM)\supseteq N.$
 Let $g\M^*$ be the expansion of the structure $g\M$ by all the relations definable in $\NN^{Eq}.$
 
 The inclusion $\dcl(gM)\supseteq N$ allows to interpret the set $N$ as well as any relation on $N,$ in $g\M^*$ using parameters in $gM^{Eq}.$
 But the relations of $g\M^*$ are  definable in $g\M^{Eq}$ using parameters in $gM^{Eq}$ (see \ref{facts}) thus we conclude that $\NN$ is definable in $gM^{Eq}$ using some parameters $A,$ or there is a morphism
 $$h:\NN\to g\M/A.$$  
Since $A$ consists of canonical parameters of sets and relations definable in $\NN^{Eq}$ the morphism $h$ is an embedding. But it is also a surjection by construction. Hence $h$ is an $\mathfrak{M}$-isomorphism.
 $\Box$

 \epk 

\bpk \label{M/N} {\bf Lemma.} {\em Suppose $\M\in \mathfrak{M}_\mathrm{fin}.$ 
Let $H\hookrightarrow \mathrm{Aut}(\M)$ be a closed subgroup. Then $H$ is a pointwise stabiliser of a subset $A\subset \M^{eq}$ (first-order imaginaries).  That is $$H=\mathrm{Aut}(\M/A)$$

$H$ is normal if and only if the restriction of $\dcl(A)$ to any finite $\M_\alpha$ (in the notation of \ref{finry}) is first-order $0$-definable.}

\medskip

{\bf Proof.} 
 The equality $H=\mathrm{Aut}(\M/A)$  follows from \ref{EqCat}(iv) and \ref{M/A}.
Since $H$ is closed in profinite topology,
$$H=\lim_{\leftarrow}H_\alpha,\ \ H_\alpha\hookrightarrow \Aut\,\M.$$
The functorial correspondence of Theorem \ref{EqCat}(iv) identifies $H_\alpha=\Aut\, \NN_\alpha$ for some finite $\NN_\alpha$ which satisfies the assumptions of \ref{M/A} and thus 
$$H_\alpha=\Aut\, \M_\alpha/A_\alpha$$
where  $A_\alpha$ are the respective imaginaries in $\M^{eq},$ which are first order since $\M_\alpha$ is finite. By functoriality of the construction 
$$A:=\lim_{\rightarrow} A_\alpha$$ has the required property.

With the above choice of $A,$  $H$ is normal iff $N$ is invariant under  $\mathrm{Aut}(\M),$ that is $A$
0-definable. $\Box$
\epk

\bpk\label{Pr1} {\bf Proposition.} {\em To every $0$-definable  $\NN$ in $\M$ (write $\NN\hookrightarrow \M$)
one associates the exact sequence
\be\label{exact} 1\to \mathrm{Aut}(\M/N)\to \mathrm{Aut}(\M)\to \mathrm{Aut}(\NN)\to 1\ee
and every exact sequence of closed subgroups has this form for some  $\NN\hookrightarrow \M.$}

{\bf Proof.} Assuming $\NN\hookrightarrow \M,$
the surjection $\mathrm{Aut}(\M)\to \mathrm{Aut}(\NN)$ is just \ref{EqCat}(ii). The kernel of the latter homomorphism is clearly $\mathrm{Aut}(\M/N)$ which is normal as noticed above.

The inverse follows from \ref{M/N}. $\Box$
\epk

\bpk {\bf Lemma.} {\em 
Let $\M\in \mathfrak{M}_{\mathrm{fin}}.$ 
Then 
$\M$ is first-order homogeneous, i.e.
 for any two sequences $a,a'$ in $\M$ 
the first-order types of $a$ and $a'$ are equal if and only if there is $\sigma\in \Aut\,\M$ such that $\sigma(a)= a'.$ }

{\bf Proof.}  Since $\M$ is finitary we have $$M=\acl(\emptyset)=\acl(a)=\acl(a').$$
The condition $\tp(a)=\tp(a')$ implies the existence of an elementary monomorphism (partial isomorphism preserving all first-order formulas)
$\sigma: a\mapsto a'.$ It is a standard fact that any elementary monomorphism can be lifted to monomorphism $\acl(a)\to \acl(a').$
$\Box$
\epk
\section{Sections and section-imaginaries}

 
\bpk\label{sections}
 {\bf Theorem.} {\em Let $\mathbf{N}, \M\in \mathfrak{M}_{\mathrm{fin}}$
be the members of the exact sequence (\ref{exact}), 
  $$\hat{h}:\mathrm{Aut}(\M)\twoheadrightarrow\mathrm{Aut}(\mathbf{N}).$$
  
Suppose  there exists
$\hat{g}:\mathrm{Aut}(\mathbf{N})\hookrightarrow \mathrm{Aut}(\M),$  a section of 
$\hat{h},$ that is 
$$\hat{h}\circ \hat{g}=\mathrm{id}_{\mathrm{Aut}(\mathbf{N})}.$$

Then there exists a  set of first-order imaginaries $A\subset \M^{eq}$ such that the  embedding of structures $h\NN\subseteq \M^{eq}$ associated with $\hat{h}$ gives rise to the 
 $\mathfrak{M}$-isomorphism \be\label{sections0} g: h\mathbf{N}\cong \M/A,\ee
 
satisfying the following two conditions:
 \be\label{sections1} M\subseteq\dcl_{\M^{eq}}(hN\cup A)\ee
 and
\be \label{sections2} 
    \dcl_{\M^{eq}}(A)\cap \dcl_{h\NN^{eq}}(hN) =\dcl_{h\NN^{eq}}(\emptyset) .\ee

Conversely,  suppose there exist $A$ and  an interpretation-isomorphism 
(\ref{sections0}) which satisfy (\ref{sections1}) and (\ref{sections2}). Then the respective homomorphism $$\hat{g}:\mathrm{Aut}(\mathbf{N})\hookrightarrow \mathrm{Aut}(\M)$$ is a section of 
$\hat{h}.$

   }
   
   {\bf Proof.} By the assumptions we are also given an interpretation 
  $$h: \NN\hookrightarrow \M$$ 
  correponding to $\hat{h},$ such that any $\sigma\in \Aut(\M)$ induces  $\hat{h}(\sigma)\in \Aut(h\NN)$ and in this way we get all automorphisms of $h\NN,$ that is $\hat{h}(\Aut(\M))=\Aut(h\NN).$ Without loss of generality we may assume that $\NN$ is a substructure of $\M^{Eq},$ that is
  $h$ is a pointwise identity embedding and $\hat{h}(\sigma)$ is the restriction of $\sigma$ to $\NN.$ Thus 
  
  \be\label{hN} \hat{h}(\Aut(\M))=\Aut(\NN).\ee 
    
   Consider the subgroup $\hat{g}(\Aut(\NN))\subseteq \Aut(\M),$ an isomorphic copy of $\Aut(\NN).$ Since $\hat{g}$ is a section of $\hat{h}$ we get $$\hat{h} (\hat{g}(\Aut(\NN)))=\Aut(\NN).$$
   
By assumptions $\hat{g}$ lifts any automorphism
$\rho\in \Aut(\NN)$ to a unique automorphism
$\hat{g}(\rho)\in \Aut(\M)$ giving the embedding $\hat{g}:\Aut(\NN)\hookrightarrow \Aut(\M).$ 
   
   
Set $$A:= \mathrm{Fix}_{\M^{Eq}}(\hat{g}(\Aut(\NN))).$$

Note that according to \ref{M/N}
$$\hat{g}(\Aut(\NN))=\Aut(\M/A)$$  
     Moreover, $\NN$ is definable in $\M^{eq}$ over $A$ since $\NN$ as a structure is  $\Aut(\M^{eq}/A)$-invariant. In other words there is a pre-morphism $$g:\NN\to \M$$
     realised by the  embedding of the universe into $\M^{eq},$ i.e. $g(x)=x$ for any $x\in N.$    
   
Now note that $g$ is an $\mathfrak{M}$-embedding since every  $\rho\in\Aut(\NN)$ lifts to an automorphism $\hat{g}(\rho)$ of $\M^{eq}/A.$

Next we note that $g$ is an $\mathfrak{M}$-surjection, that is $\dcl_{\M^{eq}/A}(N)\supseteq M,$ or equivalently  $$\dcl_{\M^{eq}}(N\cup A)\supseteq M.$$ To see the latter we remark that if $\sigma\in \Aut(\M)$ fixes $A\cup N$ point-wise then $\sigma\in \hat{g}(\Aut(\NN))$ (because $A$ is fixed) and $\sigma$ is identity on $N.$ That is $\sigma=\mathrm{id}.$  

Finally note that $$A\cap\dcl_{\NN^{eq}}(N) =\dcl_{\NN^{eq}}(\emptyset)=\dcl_{\M^{eq}}(\emptyset)\cap \dcl_{\NN^{eq}}(N), $$        
the first equality follows from the fact that the intersection consists exactly of $\Aut(\NN)$-fixed points of $\NN^{eq},$ and the second equality is the consequence of $g$ being an embedding.  $\Box$  
\epk 
 \bpk We call $A$ satisfying (\ref{sections1}) and (\ref{sections2}) of \ref{sections} a  {\bf section-imaginary}, or more precisely, the section-imaginary corresponding to the interpretation-isomorphism $g$ of (\ref{sections0}) and section $\hat{g}$ of $\hat{h}.$
 
  Note that by definition 
$$A= \mathrm{Fix}_{\M^{eq}}(\hat{g}(\Aut(\NN)))$$
and so is totally determined by the morphism $g: \NN\to \M.$ 
\epk 
\section{Grothendieck's anabelian section conjecture}
\bpk  The celebrated Grothendieck's section conjecture is formulated in terms of a smooth algebraic curve $\X$ defined over a number field $\kk,$  its \'etale fundamental group $\pi_1^{et}(\X)$ and a section of the canonical surjective homomorphism $\pi_1^{et}(\X)\twoheadrightarrow \mathrm{Gal}_\kk.$

 In our setting $\X(\bar{\kk})$ corresponds to $\NN$ in the language that has names for every point of $\kk$ (in particular, definable points in $\NN$ are exactly $\kk$-rational points).  $\pi_1^{et}(\X)$ corrsponds to $\Aut\, \M$ and $\mathrm{Gal}_\kk$ to $\Aut\, \NN.$ As explained in \cite{AZ} in Grothendieck's setting for $\M$ one can take a multisorted cover structure $\tilde{\X}^{et},$ so our $\M$ should be seen as just one layer of  $\tilde{\X}^{et}.$
  However, we believe that just one high enough layer of the cover suffices to detect the existence of a rational point on $\X.$
 
 Note that Grothendieck also assumed that $\X$ is "anabelian", more concretely, of genus $>1.$ This is the condition necessary for the correspondence \begin{center}
 (conjugacy classes of) sections $\leftrightarrow$ rational points  
\end{center}
to be bijective, see \cite{CStix}. We only consider the conjecture
   \begin{center}
 existence of sections $\leftrightarrow$ existence of rational points  
\end{center}
which makes sense for much broader class of varieties $\X$ and is known, see \cite{Stix}, Cor.102, to imply the validity of  the original Grothendieck conjecture. In particular, the latter version of the conjecture is open for curves of genus 1, which from model theory point of view are $E$-torsors, for $E$ a group structure of an elliptic curve. See \cite{Stix}, \cite{CStix} and \cite{ZGSC} for some results on this case. The latter text also presents  Grothendieck's section conjecture in a way fitted better for the setting of the next section.
\epk
\section{Elimination of section-imaginaries}
In this section we  will be careful to distinguish between $\M$ and $\M^{Eq}.$

 Our aim here is to study conditions for the existence of sections in cases which could be seen as general model-theory style analogues of finite \'etale covers of smooth algebraic varieties.  

\bpk 
We consider specific two-sorted {\bf cover structures} 
$$\M=(\mathbf{C}, \NN, \pr)$$
 where $\mathbf{C}$ and $\NN$ are substructures on universes $C$ and $N$ respectively and 
 $$\pr: C\twoheadrightarrow N$$
 is a covering map with finite fibres. We also assume that  the image $\pr(R)$ of any definable relation $R\subseteq C^n$ is already definable in  $\NN.$ 

This gives us 
an interpretation-embedding
$$i_\pr: \NN\hookrightarrow_\mathfrak{M} \M, $$
given by the identity map on $N,$ $i_\pr N=N,$ and the corresponding short exact sequence (\ref{exact}) in Proposition \ref{Pr1}. 

Assume further that:
\begin{itemize}
\item[C1.] For any $c\in C,$
$$\dcl(N\cup\{ c\})\supseteq C.$$ 
\item[C2.] There is a 0-definable group $\Gamma$ of transformations of $C$ which fixes fibres $\pr\inv(n),$ $n\in N,$ and  acts freely and transitively on each fibre. Moreover,
 $\Gamma=\Aut(\M/N)$ (the {\bf deck-group}).

\item[C3.] The first-order theory of $\M$ is categorical in uncountable cardinals.
\end{itemize}

\epk
\bpk \label{CC1}
Note that C1 brings us into the context of {\bf generalised imaginaries} of Hrushovski's paper \cite{HrGI}. In particular, the construction there of the {\bf definable groupoid} is applicable here: 
{\em there is a 0-definable $\mathcal{G}\subseteq \NN^{eq}$ acting partially on $C$ so that:
 
for each $c_1,c_2\in C$ there is a unique   $g\in \mathcal{G}$ satisfying 
$c_2=g*c_1,$ and 

for each $g\in \mathcal{G}$ there are $c_1,c_2\in C$ 
 satisfying 
$c_2=g*c_1.$}

\medskip

The deck-group  $\Gamma$ is a subgroupoid of $\mathcal{G}$ and also is the liaison group of $C$ over $N.$

\medskip

One of the essential differences in the approaches here and in \cite{HrGI} is that the main interest of the latter is in the case when $\NN$ eliminates (ordinary) imaginaries which in the context  of algebraic geometry corresponds to $\NN$ being of genus 0.

\epk
\bpk\label{Exm} {\bf Example. A curve of genus 1.} Consider a projective curve $\X_a$ over $\Q$ given by the homogeneous equation 
$$ x^3+y^3=az^3\mbox{ some }a\in \Z$$
For almost all $a$ the equation defines a smooth complex curve and we assume $a$ has this property. For some $a,$ e.g. $a=3$ the equation
 has no non-trivial solutions in $\Q,$
that is $\X_a(\Q)=\emptyset.$ According to the general theory the quibic equation defines over $\bar{\Q}$  an elliptic curve with a well-known group structure $$\mathcal{G}_a=(\X_a(\bar{\Q}),+, 0)$$ which is  determined by the choice of the point $0\in \X_a(\bar{\Q}).$ In fact, $(\X_a(\bar{\Q}),+, 0)$ is isomorphic to $(\X_a(\bar{\Q}),+', 0')$ by  
$$ x+' y:= x+y-0',$$ for any other choice of a zero $0'\in \X_a(\bar{\Q}).$  It follows that the group $\mathcal{G}_a$ is definable over $\Q.$

Now let $\NN_a$ be the structure on the universe $N_a:=\X_a(\bar{\Q})$ in the relational language $\LL_{\NN}$ which consists of all relations $R\subseteq N_a^n$
that are Zariski closed and defined over $\Q.$ It is immediate that
$\mathcal{G}_a$ is interpretable in $\NN_a$ with $N_a$ as the universe of $\mathcal{G}_a$ 
 and that $\mathcal{G}_a$ acts on $N_a$ freely and transitively, that is $(N_a, \mathcal{G}_a)$ is a  torsor.

Let $C_a=\{ c_g: g\in \mathcal{G}_a\}$ and define the torsor
$\mathbf{C}_a=(C_a, \mathcal{G}_a)$ by the action $$(h, c_g)\mapsto c_{g+h}, \ h\in \mathcal{G}_a, c_g\in C_a.$$

Choose an $n>1$ and define $$\Gamma_{a,n}=\{ \gamma\in \mathcal{G}_a: n\cdot \gamma=0\}.$$
Finally define $$\pr_n: C_a\to N_a; \ c_g\mapsto n\cdot g,$$ 
$$\M_{a,n}=(\mathbf{C}_a, \NN_a, \pr_n).$$

This satisfies C1--C3.

\epk
\bpk
{\bf Lemma.} {\em Assume a cover structure 
$(\mathbf{C}, \NN, \pr)=\M\in \mathfrak{M}_{fin}$ and satisfies C1 and C2.

Suppose there exists $e\in \dcl(\emptyset)\cap N.$ Then 
for any $a\in \pr\inv(e)$ the set
$A=\dcl(a)$ is a section-imaginary for some section $\hat{g}: \Aut\, \NN\to \Aut\, \M$ of $\hat{i}_\pr.$
}

{\bf Proof.} We may consider the interpretation $i_\pr$ as also being an interpretation of $\NN$ in $\M/A,$ call it $i_\pr^A$ in this context.
  
Condition C1 implies that $\dcl(N\cup A)$ contains all the points of universes of $\M/A$ and thus we have (\ref{sections1}) satisfied. 
$i_\pr^A$ is a surjection in category $\mathfrak{M}.$

In order to see that $i_\pr^A$ is an embedding  consider a relation $R$ on $N$ defined by a formula  in the language of $\M/A,$ that is by a formula
$\varphi(a,y),$ where $\varphi(x,y)$ is in the language of $\M,$ $a\subseteq  A.$  It follows from C2 that, for all $a'\in \pr\inv(e),$ formula $\varphi(a',y)$ defines the same relation $R$ on $N,$ that is $R$ is defined in the language of $\M.$ By our assumption $R$ is then defined in the language of $\NN.$ Thus $i_\pr^A$ is an embedding   in category $\mathfrak{M}$ and (\ref{sections2}) is satisfied. $\Box$
\epk
\bpk Given a cover structure $(\mathbf{C}, \NN, \pr),$ we say that {\bf section-imaginaries for $i_\pr$ are eliminable} if 
 any section-imaginary $A$ for $i_\pr$ is of the form $A=\dcl(A_0),$ for some $A_0\subseteq C$.
 
 Equivalently, it follows from C1 that section-imaginaries for $i_\pr$ are eliminable if and only if 
 any section-imaginary $A$ for $i_\pr$ is of the form $A=\dcl(a),$ for some $a\in C$. 
\epk
\bpk {\bf Lemma.} {\em Assume a cover structure 
$(\mathbf{C}, \NN, \pr)=\M\in \mathfrak{M}_{fin}$ and the cover is  non-trivial, i.e. $|\Aut(\M/N)|>1.$
 Suppose  there is a section $\hat{g}: \Aut\, \NN\to \Aut\, \M$ of $\hat{i}_\pr$ and
section-imaginaries for $i_\pr$ are eliminable. Then 
$$\dcl(\emptyset)\cap N\neq \emptyset.$$ 
}

{\bf Proof.} Let $A=\dcl(A_0)$ be the section-imaginary for $i_\pr,$
$A_0\subset C.$ Note that $A_0\neq \emptyset$ since by (\ref{sections1})
$C\subset \dcl(A_0\cup N)$ and  $C\nsubset \dcl(N)$  by condition C2.

Clearly, $\pr(A_0)\subset N$ and  $\pr(A_0)\subset \dcl(A_0).$ It follows from condition (\ref{sections2}) on section-imaginaris that $\pr(A_0)\subset \dcl(\emptyset).$ $\Box$
\epk

\bpk \label{Main} {\bf Theorem.} {\em Let  
$(\mathbf{C}, \NN, \pr)=\M\in \mathfrak{M}_{fin}$ be a non-trivial cover structure and assume $\M$ has elimination of section-imaginaries for $\hat{i}_\pr$ . Then there is a section $\hat{g}: \Aut\, \NN\to \Aut\, \M$ of $\hat{i}_\pr$ if and only if there is a definable point in $\NN.$
} 

{\bf Proof.} Follows from the two lemmas above. $\Box$
\epk
\bpk Consider the special case of $\M$ where $\NN$ is the structure on a projective curve over a number field $\kk.$ We will refer to this as  the {\bf geometric case.} The geometric case has the {\bf anabelian} version, that is the case of a curve of genus $> 1$ over a finitely generated field $\kk,$
 and the {\bf abelian} version, the curve of genus 1 as e.g. described in the Example \ref{Exm}. 
 
 Suppose $\M_{a,n}$ of  the Example, for some $a$ and $n,$ has a section-imaginary, equivalently (see Theorem \ref{sections}) there is a section of (\ref{exact}). Then we expect 
that $\X_a$ has a rational point, which will happen if the section-imaginary is eliminable. 
 
\epk
\bpk {\bf Problem.} {\em Formulate a model-theoretic condition sufficient for elimination of section-imaginaries and consistent with the geometric case}. 
\medskip

The analysis of the similar structure in \cite{AZ} shows that C1-C3 are satisfied in the geometric case.

A possible extension of these conditions could be e.g.:

C4. $\Gamma\cap \dcl(\emptyset)=\dcl_\Gamma(\emptyset)$ (the dcl in the group structure $(\Gamma,\cdot).$

\epk

\end{document}